\documentclass[10pt,journal]{IEEEtran}
\ifCLASSINFOpdf
\else
\fi

\usepackage{graphicx}
\usepackage{subfigure}
\usepackage{mathptmx} 
\usepackage{amsmath} 
\usepackage{amssymb}  

\usepackage{amsfonts}
\usepackage{color}
\usepackage{url}
\usepackage{mathrsfs}

\usepackage{epstopdf}

\newtheorem{theorem}{Theorem}[section]
\newtheorem{lemma}{Lemma}[section]

\newtheorem{remark}{Remark}[section]

\def\disp{\displaystyle}

\newcommand{\R}{{\mathbb R}}
\def\crr{\cr\noalign{\vskip2mm}}

\def\dref#1{(\ref{#1})}

\def\<{\langle}
\def\>{\rangle}

\begin{document}
%
\title{ Frequency Energy Multiplier Approach to Uniform \\
Exponential Stability Analysis of Semi-discrete \\ Scheme for a Schr\"{o}dinger Equation\\ under Boundary Feedback }
%
%
%

\author{Bao-Zhu~Guo and  Fu~Zheng 
\thanks{ This work was supported by the National Natural Science Foundation of China under grants no.61873260, 11871117, 12131008. (Corresponding author: Fu~Zheng)}
\thanks{Bao-Zhu Guo is with Department of Mathematics and Physics, North China Electric Power University, Beijing
102206, China, and Key Laboratory of System and Control,
        Academy of Mathematics and Systems Science,
        Academia Sinica, Beijing 100190,   Email:bzguo@iss.ac.cn}
	 \thanks{Fu Zheng is with School of Science  of Hainan University, Haikou, Hainan 570228,  E-mail: fuzheng@hainanu.edu.cn}
}

%
%

\markboth{IEEE TRANSACTIONS ON AUTOMATIC CONTROL,~Vol.~xx, No.~x, August~20xx}%
{Shell \MakeLowercase{\textit{et al.}}: Bare Demo of IEEEtran.cls for IEEE Journals}
%



\maketitle

\begin{abstract}
  In this paper, we investigate the  uniform exponential
 stability of a semi-discrete scheme for  a  Schr\"{o}dinger equation under boundary feedback stabilizing control in the natural state space $L^2(0,1)$.
 This study is significant since  a time domain energy multiplier that allows proving the exponential stability of this continuous Schr\"{o}dinger system has not yet found, thus
 leading to a major mathematical challenge to semi-discretization of the PDE, an open problem for a long time.
 Although the powerful frequency domain energy multiplier approach has been used in proving  exponential stability
 for PDEs since 1980s, its use to the \emph{uniform} exponential  stability of the semi-discrete scheme for PDEs has not been reported yet.
 The difficulty associated with the uniformity is that due to the parameter of the step size, it involves a family of   operators in different state spaces that need to be considered simultaneously.
 Based on the Huang-Pr\"{u}ss frequency domain criterion for uniform exponential stability of a family of $C_0$-semigroups in Hilbert spaces,
 we solve this problem for the first time by proving the uniform boundedness for all the resolvents
 of these operators on the imaginary axis. The proof almost exactly follows the procedure for the exponential stability of the continuous counterpart, highlighting
 the advantage of this discretization method.

\end{abstract}


%
\IEEEpeerreviewmaketitle

\section{Introduction}

\IEEEPARstart{C}ontrol  systems described by partial differential equations (PDEs) is  infinite-dimensional.
Being such, its controller such as the
observer-based feedback control is also infinite-dimensional. As a result,  the discretization finds itself in almost all implementations of PDE control. Among  many discretization
methods is the finite-difference method which becames popular due to its simplicity in principle and its appeal to  engineers.
One of the most  commonly used  discretization method is the so-called semi-discrete scheme which
keeps time continuous while discretizing the spatial variable.
It has been widely studied in literature.
The main advantage of the semi-discrete
scheme is that it results in an ordinary differential equation system,  which control researchers are most familiar with.
However,  it  has  been acknowledged for a long time that the uniform exponential stability with respect to the spatial discrete step size
cannot be guaranteed for classical semi-discrete schemes for PDEs, largely due to presence of  high frequency spurious components.
In addition, some other typical important
control properties such as uniform
observability and uniform exact controllability cannot be guaranteed either.
The reason for this loss is that the spurious modes are only weakly damped in the process of semi-discretization. A detail account can be found in    \cite{Zuazua}.  For wave equations, several remedies such as  Tichonoff regularization \cite{Glowinski},  mixed-finite elements \cite{Banks,Micu},       high frequency filtering \cite{Infante}, and non-uniform meshes \cite{Ervedoza1}, have been proposed to circumvent  this difficulty.
Among many these  remedies,  the numerical viscosity damping  introduced in \cite{Tebou1,Tebou2} is the most popular.
However,   this approach brings  a viscosity term artificially added  into the classical discrete scheme.
The coefficients of the numerical viscosity damping  vary from PDE to PDE.
Recently, a new natural semi-discrete scheme based on order reduction finite
difference method was introduced in  \cite{Liujk1} and has been applied to different systems \cite{Guobz2,Zheng}.
This approach has the critical advantages that it guarantees the uniform exponential stability. In addition, as a natural semi-discrete scheme, it allows one to prove the uniform exponential stability in a manner
parallel to its continuous PDE counterpart.
Nevertheless, all the previous papers on this scheme involved construction of Lyapunov functional which the proof heavenly relies on, both for semi-discrete scheme and the continuous counterpart.

Construction of a suitable Lyapunov functional for a PDE relies on a time domain energy multiplier, which is not always available and its construction is most often very technical.
In 1980s, a frequency domain energy multiplier approach was developed
for the exponential stability initially
for a single PDE (\cite{Liuzy}).
The approach is based on
a frequency domain characterization for exponential stability of $C_0$-semigroup in Hilbert space.
Originally developed independently in \cite{Huang} and \cite{Pruss},
the result of was proved later in \cite{Liuzys}
to be valid for uniform exponential stability of a family of $C_0$-semigroups in Hilbert spaces as well.
Uniform admissibility and observability for the finite element space semi-discretizations of abstract  Schr\"{o}dinger system and
second order infinite dimensional vibrating systems have also been developed \cite{Ervedoza2,Ervedoza3}.

In this paper, we investigate the uniform exponential stability of an order reduction semi-discrete scheme for a  Schr\"{o}dinger equation under
 boundary control by the frequency domain multiplier approach. It is significant because one cannot find a suitable time domain
 Lyapunov functional both for the continuous PDE and for its discrete scheme. This implies that
successful approaches presented in \cite{Liujk2, Liujk1,Guobz2,Zheng} cannot be applied here.
As a matter of fact, in order to apply the Lyapunov
method, the paper \cite{Liujk2} has to consider the  Schr\"{o}dinger system in the high order state space $H^{1}(0,1)$,
whereas our state space is the standard space $L^2(0,1)$. The  problem in $L^2(0,1)$ has been open for quite a long time. Fairly speaking, this paper brings a new way to the proof of the uniform exponential stability of the semi-discrete scheme for PDEs.
It is also worthy pointing out that the proofs for both continuous PDE and for the discrete counterpart are again  analogous, demonstrating the advantage of the order reduction semi-discretization approach.

We proceed as follows. In the next section,  Section  \ref{sec2}, we prove the exponential stability of the continuous PDE by the frequency domain multiplier method. Although it is the simplest PDE ever studied in the literature,
it helps in constructing a frequency domain multiplier for its semi-discrete counterpart.
In  Section \ref{sec3}, we design a  semi-discretized scheme and obtain a family of finite-dimensional systems.
In Section \ref{sec4}, the uniform exponential stability is developed by the frequency domain multiplier approach.
We introduce the shadow element to  help  understand the numerical approximating scheme, which plays an important role in the proof of uniform stability.  Some concluding remarks are included in Section \ref{sec5}.

\section{Stability of Schr\"{o}dinger system  via frequency domain multiplier}\label{sec2}
Consider the  following Schr\"{o}dinger equation under  boundary control:
\begin{equation}\label{2.p1}
\left\{\begin{array}{l}
w_{t}(x,t)=-iw_{xx}(x,t),~t>0,~x\in(0,1), \crr
w(0,t)=0,~t\geq0,\crr
w_{x}(1,t)=u(t),~k>0,~t\geq0,\crr
y(t)=w(1,t),~t\ge 0,\crr
w(x,0)=w^{0}(x), ~x\in[0,1],
\end{array}\right.
\end{equation}
where $u(\cdot)$ is the control, $y(\cdot)$ is the measured
output and  $w_0(\cdot)$ is the initial state.
Under the proportional feedback control:
 \begin{equation}\label{2.1*}
 u(t)=-ki y(t), ~k>0,
\end{equation}
the closed-loop system of \dref{2.p1}  becomes
\begin{equation}\label{2.1}
\left\{\begin{array}{l}
w_{t}(x,t)=-iw_{xx}(x,t),~t>0,~x\in(0,1), \crr
w(0,t)=0,~t\geq0,\crr
w_{x}(1,t)=-ki w(1,t),~k>0,~t\geq0,\crr
w(x,0)=w^{0}(x), ~x\in[0,1],
\end{array}\right.
\end{equation}
 We consider system \dref{2.1}
in the natural   state space $L^2(0,1)$. Define the  system operator of  \dref{2.1} as follows:
\begin{equation}\label{op1}
\left\{\begin{array}{l}
Af=-if'', \forall f\in D(A),\crr
\disp D(A)=\{f\in L^2(0,1)|f\in H^{2}(0,1),\crr
\disp~~~~~f(0)=0,~f'(1)=-ikf(1)\}.
\end{array}\right.
\end{equation}
Then, \dref{2.1} can be written as an evolution equation in $L^2(0,1)$:
\begin{equation}\label{add1}
\left\{\begin{array}{l} \dot{w}(\cdot,t)=Aw(\cdot,t),\crr
w(x,0)=w_0(x).
\end{array}\right.
\end{equation}
It is seen that
\begin{equation}\label{dis1}
 \textrm{Re}\langle Af,f\rangle_{L^2(0,1)}
  =\textrm{Re}\int_{0}^{1}-ikf''(x)\overline{f(x)}dx=-k|f(1)|^{2},
\end{equation}
which implies  that $A$ is dissipative. In addition, the operator $A$ is invertible and
\begin{equation}\label{in1}
\begin{array}{l}
A^{-1}f(x)=\disp \frac{-kx\int_{0}^{1}xf(x)dx}{1+ki}\crr
\disp ~~~~~~~~~~~~~~-i\int_{x}^{1}(x-\tau)f(\tau)d\tau-i\int_{0}^{1}xf(x)dx,
\end{array}
\end{equation}
which is bounded in $L^2(0,1)$. As a result,  $A$ generates a $C_0$-semigroup of contractions on $L^2(0,1)$ by the Lumer-Phillips theorem (\cite[Theorem 3.8.4]{Tucsnak}) and since $A^{-1}$ is compact, the spectrum of $A$ consists of
isolated eigenvalues only.

Furthermore, define the
system energy for \dref{2.1} as
\begin{equation}\label{2.1**}
 E(t)=\frac{1}{2}\int_0^1|w(x,t)|^2dt,
\end{equation}
which is non-increasing as a consequence of \dref{dis1}:
\begin{equation}\label{2.1***}
 \dot{E}(t)=-k|w(1,t)|^2.
\end{equation}

We point out that a different version of \dref{2.1}:
\begin{equation}\label{liu}
\left\{\begin{array}{l}
w_{t}(x,t)=-iw_{xx}(x,t),~t>0,~x\in(0,1), \crr
w(0,t)=0,~t\geq0,\crr
w_{x}(1,t)=-k w_t(1,t),~k>0,~t\geq0,\crr
w(x,0)=w^{0}(x), ~x\in[0,1],
\end{array}\right.
\end{equation}
was investigated in \cite{Liujk2}, for which one can find a time  domain energy multiplier,  and
a Lyapunov functional was then constructed to both system \dref{liu} and its semi-discrete
counterpart.  However, system \dref{2.1} is  a rather unusual system for which a time domain energy multiplier has been not found yet. A first exponential stability result
of system \dref{2.1} was  proved by the Riesz basis approach in \cite{Kristic2}. Although the Riesz basis
is powerful and the result obtained is much deeper than the result obtained from the multiplier method; for instance the
spectrum-determined growth condition is usually a consequence of the Riesz basis approach yet this is usually not the case with the multiplier method. Unfortunately,
the Riesz basis is extremely difficult,  at least at the moment, to be applied for the uniformly exponential stability of semi-discrete model for \dref{2.1} developed in this paper.

In this paper, we use an alternative
powerful method called the frequency energy multiplier method, which has been developed for continuous
PDEs over the last three decades  (\cite{Liuzy}). In stability analysis, we can almost give one-to-one correspondence
from continuous system to its discrete counterpart using this method. Our approach  is so powerful that can be applied to
other PDEs as well.
For notation simplicity, hereafter, we omit  without confusion the obvious dependency in time and spatial domains. The
$\mathbb{C}^n$ denotes the $n$-dimensional complex Euclidean space;  the $\mathbb{N}^+$ stands for the
set of the positive integer numbers; and $\mathbb{R}$ the set of  real numbers.

Since   $A$ generates a $C_0$-semigroup of contractions on $L^2(0,1)$, a well-known result of Hung-Pr\"{u}ss theorem \cite{Huang,Pruss} states that the  $C_{0}$-semigroup
generated by $A$ is exponentially stable if and only if
it possesses the following two properties:
\begin{enumerate}
    \item   Every imaginary number belongs to the resolvent set of  $A$, that is, $i\R\subset\rho(A)$.
    \item    The inverse operator of $i\omega-A$ is uniformly bounded for all imaginary numbers, that is,
\begin{equation}\label{1.3}
  \sup_{\omega\in \R}\|(i\omega-A)^{-1}\|<\infty.
\end{equation}
\end{enumerate}

The property $i\R\subset\rho(A)$ is stated in the following Lemma \ref{pro2.1}.
\begin{lemma}\label{pro2.1} Let $A$ be defined by \dref{op1}. Then,
$i\R\subset \rho(A)$.
\end{lemma}

{\bf Proof.} If there exist $\beta\in\mathbb{R}, \beta\neq 0$ and a nonzero $f\in D(A)$ such that $i\beta f=Af$, then
\begin{equation}\label{1.8*}
\left\{\begin{array}{l}
  i\beta f(x)=-if''(x),\crr\disp
  f'(1)=-ki f(1),~ f(0)=0.
  \end{array}\right.
\end{equation}
Take the inner product with $f(\cdot)$ over $[0,1]$ on both sides of the first equation of \dref{1.8*} to obtain
\begin{equation}\label{1.8}
  i\beta \|f\|^2=-k |f(1)|^2+i\int_0^1|f'(x)|^2dx,
\end{equation}
which gives $f(1)=0$ and hence $f'(1)=0$. This shows that \dref{1.8*} has only zero solution, a contradiction.
 \hfill\rule{2 mm}{3mm}

 \vspace{0.2cm}

\begin{theorem}\label{Th.A1}
 Let $A$ be defined by \dref{op1}. Then, \dref{1.3} holds true. As a consequence, the $C_0$-semigroup
 $e^{At}$ generated by $A$ is exponentially stable in $L^2(0,1)$.
 \end{theorem}

{\bf Proof.}   We prove   by assuming contrary of  (\ref{1.3}) that  there exit a sequence
$\omega_{n}\rightarrow\infty$,  $f_{n}\in D(A)$, $\|f_{n}\|=1$  that
\begin{equation*}
  \lim\limits_{n\rightarrow\infty}\|(i\omega_{n}-A)f_{n}\|=0,
\end{equation*}
i.e.,
\begin{equation}\label{1.4}
  i\omega_{n}f_{n}+if_{n}''\rightarrow 0 \hbox{ in } L^2(0,1).
\end{equation}
Since
\begin{equation}\label{1.5}
  \textrm{Re}\langle(i\omega_{n}-A)f_{n},f_{n}\rangle_{L^2(0,1)}=\textrm{Re}\langle-Af_{n},f_{n}\rangle
  =k|f_{n}(1)|^{2}\rightarrow0,
\end{equation}
by the  boundary condition $f'_{n}(1)=-ikf_{n}(1)$, it  gives
\begin{equation}\label{1.6}
  f'_{n}(1)\rightarrow0.
\end{equation}
From (\ref{1.4}) and $\|f_{n}\|=1$, it follows that $\frac{f''_{n}(\cdot)}{\omega_{n}}$ is bounded in $L^2(0,1)$. By
\begin{equation*}
  |f'_{n}(x)-f'_{n}(1)|=\left|\int_{1}^{x}f_n''(s)ds\right|\leq\|f_n''\|,
\end{equation*}
  it follows from (\ref{1.6}) and $\omega_{n}\rightarrow\infty$  that
\begin{equation}\label{1.7}
  \frac{f'_{n}(\cdot)}{\omega_{n}}~\textrm{is bounded in}~L^2(0,1).
\end{equation}
Since
\begin{equation*}
\begin{array}{l}
\disp {\rm Re} \left\langle\omega_{n}f_{n}+f''_{n},\frac{xf'_{n}}{\omega_{n}}\right\rangle_{L^2(0,1)}
=\frac{|f_{n}(1)|^{2}}{2}-\frac{1}{2}\int_{0}^{1}|f_{n}(x)|^{2}dx\crr\disp
~~~~+\frac{1}{2\omega_{n}}|f'_{n}(1)|^{2}-\frac{1}{2\omega_{n}}\int_{0}^{1}|f'_{n}(x)|^{2}dx,\end{array}
\end{equation*}
and
$$
\left\langle\omega_{n}f_{n}+f''_{n},\frac{xf'_{n}}{\omega_{n}}\right\rangle_{L^2(0,1)}\to 0,
$$
we  have by (\ref{1.5}) and (\ref{1.6}) that
\begin{equation}\label{1.7+}\int_{0}^{1}|f_{n}(x)|^{2}dx+\frac{1}{\omega_{n}}\int_{0}^{1}|f'_{n}(x)|^{2}dx\rightarrow0,
\end{equation}
which shows that when $\omega_n>0$,  $\|f_{n}\|^{2}\rightarrow0$, which is  a contradiction to $\|f_{n}\|=1$. On the other hand, since
from \dref{1.4} and $\omega_n\to \infty$, we have
\begin{equation}\label{A4}
\begin{array}{l}
\disp
\quad\int_{0}^{1}\left|f_{n}(x)+\frac{f''_{n}(x)}{\omega_{n}}\right|^{2}dx \crr
\disp
=\int_{0}^{1}\left(f_{n}(x)+\frac{f''_{n}(x)}{\omega_{n}}\right)\left(\overline{f_{n}(x)}+\frac{\overline{f''_{n}(x)}}{\omega_{n}}\right)dx\crr
\disp =\int_{0}^{1}\left[|f_{n}(x)|^{2}+\frac{|f''_{n}(x)|^{2}}{\omega_{n}^{2}}\right]dx\crr
\disp~~~
+\frac{1}{\omega_{n}}\int_{0}^{1}[f_{n}(x)\overline{f''_{n}(x)}+\overline{f_{n}(x)}f_{n}''(x)]dx\crr
\disp =\int_{0}^{1}\left[|f_{n}(x)|^{2}+\frac{|f''_{n}(x)|^{2}}{\omega_{n}^{2}}\right]dx\crr
\disp~~~
+\frac{1}{\omega_{n}}[f_{n}(x)\overline{f'_{n}(x)}+\overline{f_{n}(x)}f_{n}'(x)]_{0}^{1}\crr
\disp~~~-\frac{2}{\omega_{n}}\int_{0}^{1}|f'_{n}(x)|^{2}dx\to 0,
\end{array}
\end{equation}
Substitute  $f'_{n}(1)=-ikf_{n}(1)$ and $f_{n}(0)=0$ into  \dref{A4},    and use   (\ref{1.5})-(\ref{1.6}) to obtain
\begin{equation}\label{1.7++}
\int_{0}^{1}|f_{n}(x)|^{2}dx+\int_{0}^{1}\frac{|f''_{n}(x)|^{2}}{\omega_{n}^{2}}dx
-\frac{2}{\omega_{n}}\int_{0}^{1}|f'_{n}(x)|^{2}dx\rightarrow0.
\end{equation}
which shows that when  $\omega_{n}<0$,   $\|f_{n}\|^{2}\rightarrow0$, which is  also a contradiction.
  \hfill\rule{2 mm}{3mm}

  \section{Semi-discrete  scheme of Schr\"{o}dinger equation}\label{sec3}

  In this section we apply the order reduction method to derive a semi-discrete scheme for \dref{2.1}.
  To this purpose, we introduce an intermediate variable $v(x,t)=w_{x}(x,t)$ to reduce the order of the spacial derivative of
    (\ref{2.1}).  In this way, the Schr\"{o}dinger equation (\ref{2.1}) can be rewritten as the following equivalent form:
\begin{equation}\label{3.1}
\left\{\begin{array}{l}
w_{t}(x,t)+iv_{x}(x,t)=0,  \\
v(x,t)=w_{x}(x,t),\\
w(0,t)=0,\\
v(1,t)=-kiw(1,t),\\
w(x,0)=w^{0}(x).
\end{array}\right.\
\end{equation}

The semi-discretization process is similar to \cite{Liujk2}. For the sake of completeness, we sketch briefly the
 process.  For fixed \(N\in \mathbb {N}^{+}\), consider an equidistant partition of interval \([0,1]\):
$$
0=x_{0}<x_{1}<\cdots<x_{j}=jh<\cdots<x_{N+1}=1,
$$
where \( h=\frac{1}{N+1}\) is the mesh size. Denote the sequence $ \{u_j\}_{0}^{N+1}$ by $ \{u_j\}_j$ and  introduce respectively
 the average  operator  and  the first-order
 finite difference operator as
\begin{equation}\label{oo1}
u_{j+\frac{1}{2}} =\frac{u_j +u_{j+1} }{2},\quad
\delta_x u_{j+\frac{1}{2}}=\frac{u_{j+1}-u_j}{h}.
\end{equation}
 For the solutions  $v(x,t)$ and $w(x,t)$ of  \eqref{3.1}, let $ \{V_j(t)\}_j$ and $\{W_j(t)\}_j$
 be grid functions at grids $ \{x_j\}_j $, satisfying
$$
V_j(t)=v(x_j,t),\quad W_j(t)=w(x_j,t),\quad 0\leq j\leq {N+1}.
$$
The first equation of system \eqref{3.1} holds at $ (x_{j+\frac{1}{2}},t) $, i.e.,
$$
w'(x_{j+\frac{1}{2}},t)+iv_x(x_{j+\frac{1}{2}},t)=0,
$$
where $x_{j+\frac12}=(j+\frac12)h$. Hereafter the prime ``$\prime $''
 represents the derivative with respect to time $t$. Replace the
 differential operator $\partial_x$ with  difference operator $\delta_x$ to get
\begin{align}\label{q4}
	 W_{j+\frac{1}{2}}'(t)+i\delta_xV_{j+\frac{1}{2}}(t)=\mathcal{O}(h^2).
\end{align}
 Similarly, for the second equation of system \eqref{3.1}, it has
\begin{align}\label{q5}
V_{j+\frac{1}{2}}(t)-\delta_xW_{j+\frac{1}{2}} (t)=\mathcal{O}(h^2).
\end{align}
 By dropping the infinitesimal terms in \eqref{q4} and \eqref{q5},
  and replacing $ W_j(t) $ and $ V_j(t) $ by $ w_j(t) $ and $ v_j(t) $, respectively,
  we arrive at a semi-discretized finite difference scheme of system \eqref{3.1} as follows:
\begin{equation}\label{2.2}
\left\{\begin{array}{l}
w'_{j+\frac{1}{2}}(t)+i\delta_{x}v_{j+\frac{1}{2}}(t)=0,~0\leq j\leq N,\\
v_{j+\frac{1}{2}}(t)=\delta_{x}w_{j+\frac{1}{2}}(t),~0\leq j\leq N,\\
v_{N+1}(t)=-kiw_{N+1}(t),~t\geq0\\
w_{0}(t)=0,\\
w_{j}(0) =w^{0}_{j},~0\leq j\leq N+1,
\end{array}
\right.
\end{equation}
where $v_{j}(t)$ and $w_{j}(t)$  are grid functions at grids $x_{j}~(0\leq j\leq N+1)$, and $w^{0}_{j}$ is the approximation of the initial value  $w^{0}(x_{j})$.

\begin{remark}\label{rem2.1}
The semi-discretized system (\ref{2.2}) is a family of differentiation-algebra systems, which is  called
singular systems for which there are huge amount of references related to them. See for instance \cite{Dai}, \cite{Kao}, \cite{Wang} and the references therein.
\end{remark}

Now,  we  eliminate the intermediate variables $v_{j}(t)$ from (\ref{2.2}).
 To  this purpose, let $$W_{h}(t)=(w_{1}(t),w_{2}(t),\cdots, w_{N+1}(t))^\top$$
be unknown variable of (\ref{2.2}) and $$V_{h}(t)=(v_{0}(t),v_{1}(t),\cdots, v_{N}(t))^\top$$ the auxiliary variable. We write (\ref{2.2})  into   vectorial form:
\begin{equation}\label{2.3}
\left\{\begin{array}{l}
\disp
D_{h}W'_{h}(t)=\disp -iM_{h}V_{h}(t)-\left(\begin{array}{l}~~~~~~~~0\\
~~~~~~~~\vdots\\
~~~~~~~~0\\
kh^{-1}w_{N+1}(t)
\end{array}
\right), \crr\disp
D_{h}^{\top}V_{h}(t)
= -M_{h}^{\top}W_{h}(t)+\left(\begin{array}{l}~~~~~~~~0\\
~~~~~~~~\vdots\\
~~~~~~~~0\\
i2^{-1}kw_{N+1}(t)
\end{array}
\right),\crr
~~W_{h}(0)\disp =(w_{0}^{0}, w_{1}^{0},\cdots,w^{0}_{N})^\top,
\end{array}\right.
\end{equation}
where the matrices  $D_{h}$ and $M_{h}$ are given by
\begin{equation}\label{dm}
\begin{array}{l}
\disp
D_{h}=\frac{1}{2}\left(\begin{matrix}
1&&&&\\
1&1&&&\\
&\ddots&\ddots&&\\
&&1&1&\\
&&&1&1
\end{matrix}\right)_{(N+1)\times(N+1)}, \crr\disp
M_{h}=\frac{1}{h}\left(\begin{matrix}
-1&1&&&\\
&-1&1&&\\
&&\ddots&\ddots&\\
&&&-1&1\\
&&&&-1
\end{matrix}\right)_{(N+1)\times(N+1)}.
\end{array}
\end{equation}
Obviously, both $D_h$ and $M_h$ are invertible. The differential algebraic system (\ref{2.2}) or (\ref{2.3}) can be written as
  an evolution equation in $\mathbb{C}^{N+1}$:
\begin{equation}\label{2.16}
\left\{\begin{array}{l}
\disp W'_{h}(t)=\mathcal{A}_{h}W_{h}(t),~~W_{h}(t)\in\mathbb{Y}_{h}=\mathbb{C}^{N+1},\crr
W_{h}(0)=(w_{1}^{0}, w_{2}^{0},\cdots,w^{0}_{N+1})^\top\in \mathbb{Y}_{h},
\end{array}\right.
\end{equation}
 where $\mathcal{A}_{h}$ is defined by
  \begin{equation}\label{A}
\begin{array}{l}
\disp  \quad \mathcal{A}_{h}Y_{h}\crr
= D_{h}^{-1}\left[iM_{h}\left(D_{h}^{\top}\right)^{-1}\left(M_{h}^{\top}Y_{h}-(0,\cdots,0,2^{-1}iky_{N+1})^{\top}\right)\right.\crr
\left.~~~~~-D_h^{-1}(0,\cdots,0,kh^{-1}y_{N+1})^{\top}\right], \crr
~~~~~~~~~~~~~~~~\forall Y_h=(y_1,y_2,\cdots,y_{N+1})^\top \in \mathbb{C}^{N+1}.\end{array}
\end{equation}

System (\ref{2.16})   is naturally discussed in the state space $\mathbb{C}^{N+1}$.
 To relate $\mathbb{C}^{N+1}$ in \dref{2.16} with the step size, we write $\mathbb{Y}_{h}=\mathbb{C}^{N+1}$ and
 define a new   inner product for $\mathbb{Y}_{h}$:
\begin{align}\left\langle Y_{h},\widetilde{Y}_{h}\right\rangle_{\mathbb{Y}_{h}}
=h\left\langle D_{h}Y_{h}, D_{h}\widetilde{Y}_{h}\right\rangle, \forall Y_{h}, \widetilde{Y}_{h}\in \mathbb{Y}_{h},
\nonumber\end{align}
where  $\langle\cdot, \cdot\rangle$ is the standard  inner product of $\mathbb{C}^{N+1}$.
 For $Y_{h}=(y_{1},y_{2},\cdots,y_{N+1})^{\top}\in \mathbb{Y}_{h}$, we choose
 the vector
 \begin{equation}\label{zh1}
 \begin{array}{l}
\disp Z_{h}=(z_{0},z_{1},\cdots,z_{N})^{\top}\in \mathbb{C}^{N+1} \hbox{  satisfying } \crr D_{h}^{\top}Z_{h}=-M_{h}^{\top}Y_{h}+(0,\cdots,0,2^{-1}iky_{N+1})^{\top}.\end{array}
 \end{equation}
   We call $Z_{h}$ the shadow element of $Y_{h}$, which can simplify significantly the notation
   in the later proofs.

 The classical semi-discrete scheme is  similar with \dref{2.16} where the average operator $D_h=I_{N+1}$, i.e.,
 \begin{equation}\label{2.16*}
\left\{\begin{array}{l}
\disp W'_{h}(t)=\mathcal{\hat{A}}_{h}W_{h}(t),~~W_{h}(t)\in\mathbb{Y}_{h}=\mathbb{C}^{N+1},\crr
W_{h}(0)=(w_{1}^{0}, w_{2}^{0},\cdots,w^{0}_{N+1})^\top\in \mathbb{Y}_{h},
\end{array}\right.
\end{equation}
 in which the  $\mathcal{\hat{A}}_{h}$ is defined by
  \begin{equation}\label{A*}
  \begin{array}{l}
\disp  \quad \mathcal{\hat{A}}_{h}Y_{h}=iM_{h}\left(M_{h}^{\top}Y_{h}-(0,\cdots,0,2^{-1}iky_{N+1})^{\top}\right)\crr
~~~~~~~~~~~~~-(0,\cdots,0,kh^{-1}y_{N+1})^{\top}.\end{array}
\end{equation}

   At the end of this section we explain the significance of the discrete scheme \dref{2.16}. We plot two figures
  in Figures 1 and   2, respectively.  Figure 1 depicts the maximal   real parts of the eigenvalues of the classical
  semi-discrete scheme \dref{2.16*} with step size $h$, from which we see that the real parts of the eigenvalues approach zero.
 Figure 2  depicts the maximal real parts of the eigenvalues of the order reduction semi-discrete scheme \dref{2.16} with the same step size,  from which we see that the real parts of the eigenvalues approach
  a negative number. In both figures, we take $k=1$.

\begin{figure}[htbp]

{
		
			\centering
			\includegraphics[width=2.6in]{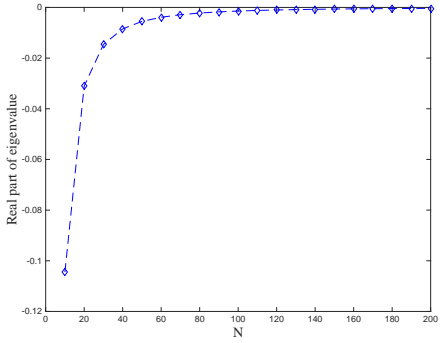}
	\centering
		}%
\caption{Maximal real parts  of eigenvalues of  the semi-discrete scheme by  classical method \dref{2.16*}}
\end{figure}\label{Fig1}

\begin{figure}[htbp]
{
				\centering
			\includegraphics[width=2.6in]{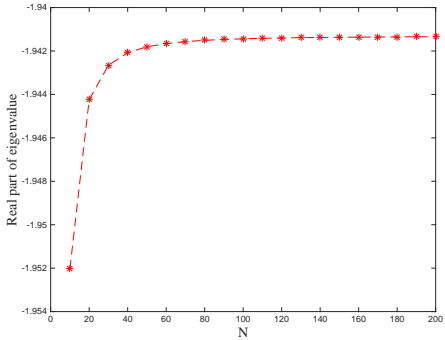}
	\centering
	}%

\caption{Maximal real parts  of eigenvalues of  the semi-discrete scheme by order reduction method \dref{2.16}}
\end{figure}\label{Fig2}

\section{Proof of Uniform exponential stability}\label{sec4}

This section is devoted to the proof of the uniform exponential stability of \dref{2.16}. To begin
with, we first show that  $\mathcal{A}_{h}$ is dissipative for every step size $h$.
\begin{lemma}\label{pro3.1+}
For the matrix $\mathcal{A}_{h}$ defined by (\ref{A}), there holds
\begin{equation}\label{3.5}
{\rm Re}\left\langle \mathcal{A}_{h}Y_{h}, Y_{h}\right\rangle_{\mathbb{Y}_{h}}
=-k|y_{N+1}|^{2},~\forall~Y_{h}\in \mathbb{Y}_{h},
\end{equation}
which implies that $\mathcal{A}_{h}$ is dissipative  for every  $h\in (0,1)$.
\end{lemma}

{\bf Proof.} For $Y_{h}=(y_{1},y_{2},\cdots,y_{N+1})\in \mathbb{Y}_{h}$, let
 $Z_{h}=(z_{0},z_{1},\cdots,z_{N})$ be the shadow element of $Y_{h}$:
\begin{equation}\label{kk1}
  \left\{\begin{array}{l}
D_{h}^{\top}Z_{h}=-M_{h}^{\top}Y_{h}+(0,\cdots,0,2^{-1}iky_{N+1})^{\top},\crr
\mathcal{A}_{h}Y_{h}=D_{h}^{-1}[-iM_{h}Z_{h}+(0,\cdots,0,kh^{-1}y_{N+1})^{\top}].
\end{array}\right.
\end{equation}
 Set $y_{0}:=0$ and $z_{N+1}:=-iky_{N+1}$ and introduce
  $\widetilde{Y}_{h}=(\widetilde{y}_{1},\widetilde{y}_{2},\cdots,\widetilde{y}_{N+1})\in \mathbb{Y}_{h}$ such that $\mathcal{A}_{h}Y_{h}=\widetilde{Y}_{h}$.   Then,
\begin{equation}\label{3.6}
 \disp D_{h}^{\top}Z_{h}+(0,\cdots,0,2^{-1}z_{N+1})^{\top}=-M_{h}^{\top}Y_{h},
\end{equation}
which is equivalent to
\begin{equation}\label{3.61}
  z_{j+\frac{1}{2}}=\delta_{x}y_{j+\frac{1}{2}},~j=0,1,\cdots,N,
  \end{equation}
  and
\begin{equation}\label{3.62}
  D_{h}\widetilde{Y}_{h}=-iM_{h}Z_{h}-(0,\cdots,0,ih^{-1}z_{N+1})^{\top},
\end{equation}
which is equivalent to
\begin{equation}\label{3.63}
\widetilde{y}_{j+\frac{1}{2}}=
-i\delta_{x}z_{j+\frac{1}{2}},~j=0,1,\cdots,N,
\end{equation}
where in all \dref{3.6} to \dref{3.63}, it was  assumed that  $\widetilde{y}_{0}=0$.
Take  the inner product between  $\mathcal{A}_{h}Y_{h}$ and   $Y_{h}$ in  $\mathbb{Y}_{h}$ by
taking \dref{3.61} and  \dref{3.63} into account  to obtain
\begin{equation}\label{gg1}
\begin{array}{l}
\disp\quad \textrm{Re}\left\langle \mathcal{A}_{h}Y_{h}, Y_{h}\right\rangle_{\mathbb{Y}_{h}}=\textrm{Re}\left\langle \widetilde{Y}_{h}, Y_{h}\right\rangle_{\mathbb{Y}_{h}}\crr=\disp
\frac{h}{2}\left\langle D_{h} \widetilde{Y}_{h}, D_{h} Y_{h}\right\rangle+\frac{h}{2}\left\langle D_{h}Y_{h} , D_{h}\widetilde{Y}_{h}\right\rangle\crr
=\disp \frac{h}{2}\sum_{j=0}^{N}\widetilde{y}_{j+\frac{1}{2}}\overline{y}_{j+\frac{1}{2}}+
\frac{h}{2}\sum_{j=0}^{N}y_{j+\frac{1}{2}}\overline{\widetilde{y}}_{j+\frac{1}{2}},~(\textrm{using } \dref{3.63})\crr
=\disp  -\frac{hi}{2}\sum_{j=0}^{N}\delta_{x}z_{j+\frac{1}{2}}\overline{y}_{j+\frac{1}{2}}
+\frac{hi}{2}\sum_{j=0}^{N}y_{j+\frac{1}{2}}\delta_{x}\overline{z}_{j+\frac{1}{2}}\crr
=\disp-\frac{hi}{2}\sum_{j=0}^{N}\left[\delta_{x}z_{j+\frac{1}{2}}\overline{y}_{j+\frac{1}{2}}(t)
+z_{j+\frac{1}{2}}\delta_{x}\overline{y}_{j+\frac{1}{2}}\right]\crr
\disp
~~+\frac{hi}{2}\sum_{j=0}^{N}\left[y_{j+\frac{1}{2}}\delta_{x}\overline{z}_{j+\frac{1}{2}}+
\delta_{x}y_{j+\frac{1}{2}}\overline{z}_{j+\frac{1}{2}}\right].~
(\textrm{using } \dref{3.61})
\end{array}
\end{equation}
A simple calculation shows that
\begin{eqnarray}\label{3.4}
\disp &&\disp -\frac{hi}{2}\sum_{j=0}^{N}\left[\delta_{x}z_{j+\frac{1}{2}}\overline{y}_{j+\frac{1}{2}}+z_{j+\frac{1}{2}}\delta_{x}\overline{y}_{j+\frac{1}{2}}\right]
\crr
&&\disp\quad+\frac{hi}{2}\sum_{j=0}^{N}\left[y_{j+\frac{1}{2}}\delta_{x}\overline{z}_{j+\frac{1}{2}}+
\delta_{x}y_{j+\frac{1}{2}}\overline{z}_{j+\frac{1}{2}}\right]\crr
&&\disp=-\frac{i}{4}\sum_{j=0}^{N}\left[(z_{j+1}-z_{j})(\overline{y}_{j+1}+\overline{y}_{j})
+(z_{j+1}+z_{j})(\overline{y}_{j+1}-\overline{y}_{j})\right]\crr
&&\disp+\frac{i}{4}\sum_{j=0}^{N}\left[(y_{j+1}+y_{j})(\overline{z}_{j+1}-\overline{z}_{j})
+(y_{j+1}-y_{j})(\overline{z}_{j+1}(t)+\overline{z}_{j})\right]\crr
&&\disp=-\frac{i}{2}\sum_{j=0}^{N}[z_{j+1}\overline{y}_{j+1}-z_{j}\overline{y}_{j}]
+\frac{i}{2}\sum_{j=0}^{N}[y_{j+1}\overline{z}_{j+1}-y_{j}\overline{z}_{j}]\crr
&&\disp=\frac{i}{2}[z_{0}\overline{y}_{0}-z_{N+1}\overline{y}_{N+1}]
+\frac{i}{2}[y_{N+1}\overline{z}_{N+1}-y_{0}\overline{z}_{0}]\crr
&&\disp=-k|y_{N+1}|^{2}. ~~(\hbox{ using } -iky_{N+1}=z_{N+1}\hbox{ and } y_{0}=0)
\end{eqnarray}
The (\ref{gg1}) and (\ref{3.4}) leads to  (\ref{3.5}).
 \hfill\rule{2 mm}{3mm}

   \vspace{0.2cm}

Define the energy of   (\ref{2.16})  as
\begin{equation}\label{4.1}
E_{h}(t)=\frac{h}{2}\sum_{j=0}^{N}\left|w_{j+\frac{1}{2}}(t)\right|^{2}=\frac{1}{2}\left\langle W_{h}(t), W_{h}(t)\right\rangle_{\mathbb{Y}_{h}}.
\end{equation}
which is the discretization of the continuous energy \dref{2.1**}.
The following Lemma \ref{pro3.1} is the discrete counterpart of \dref{2.1***}, which is a consequence of \dref{3.5}.

\begin{lemma}\label{pro3.1}
The $E_{h}(t)$  defined by \dref{4.1} satisfies
\begin{equation}\label{3.2}
\dot{E}_{h}(t)=-k|w_{N+1}(t)|^{2}.
\end{equation}
\end{lemma}

 The dissipativity of $\mathcal{A}_{h}$  implies that the spectral set  $\sigma(\mathcal{A}_{h})$ of $\mathcal{A}_{h}$
 is contained in the closed left half-plane of the complex plane $\mathbb{C}$. Actually,   we have more stronger result.
 Precisely, for any $0<h<1$, the spectral set  $\sigma(\mathcal{A}_{h})$ of $\mathcal{A}_{h}$  is contained in the open left half-plane of $\mathbb{C}$. This is   the following Lemma \ref{pro3.2}.
\begin{lemma}\label{pro3.2}
For every $h\in (0,1)$, $i\mathbb{R}\subset \rho(\mathcal{A}_{h})$.
\end{lemma}

 {\bf Proof.}~If there exist $\beta\in\mathbb{R}$ and nonzero $Y_{h}\in \mathbb{Y}_{h}$ such that $i\beta Y_{h}=\mathcal{A}_{h}Y_{h}$, then  it follows from (\ref{3.5}) that
\begin{align}
0=\textrm{Re}\left\langle i\beta Y_{h},~Y_{h}\right\rangle_{\mathbb{Y}_{h}}=\textrm{Re}\left\langle \mathcal{A}_{h} Y_{h},~Y_{h}\right\rangle_{\mathbb{Y}_{h}}=-k|y_{N+1}|^{2}\label{3.8}.
\end{align}
Replacing $\widetilde{Y}_{h}$ by $i\beta Y_{h}$ in   (\ref{3.62}),  we obtain
\begin{equation}\label{3.9}
\left\{\begin{array}{l}
\beta y_{j+\frac{1}{2}}+\delta_{x}z_{j+\frac{1}{2}}=0,~~~~~0\leq j\leq N,\crr\disp
z_{j+\frac{1}{2}}-\delta_{x}y_{j+\frac{1}{2}}=0,~~~~~0\leq j\leq N,
\end{array}\right.
\end{equation}
where  $Z_{h}$ is the shadow element of $Y_{h}$ defined in \dref{kk1}, $y_{0}=0$ and $z_{N+1}:=-iky_{N+1}$.
 Hence $z_{N+1}=y_{N+1}=0$ from (\ref{3.8}). Setting $j=N$ in (\ref{3.9})  yields
$$
\beta h y_{N}=2z_{N},~z_{N}=-\frac{2}{h}y_{N}.
$$
It follows that $y_{N}=z_{N}=0$ whenever  $\beta h^{2}+4$ is nonzero.
Under the condition $\beta h^{2}+4\neq0$,
suppose $z_{j+1}=y_{j+1}=0$ and solve  (\ref{3.9}) to arrive at  $z_{j}=y_{j}=0$. This gives   $Y_{h}=0$ by induction,  which is a contradiction.
On the other hand, whenever  $\beta h^{2}+4=0$, it follows from (\ref{3.9})  that
\begin{equation}\label{3.111}
\left\{\begin{array}{l}
\disp \frac{1}{h}(y_{j+1}+y_{j})=\frac{1}{2}(z_{j+1}-z_{j}),~j=0,1,\cdots,N,\crr\disp
\frac{1}{h}(y_{j+1}-y_{j})=\frac{1}{2}(z_{j+1}+z_{j}),~j=0,1,\cdots,N,
\end{array}\right.
\end{equation}
which implies that
\begin{equation}\label{3.113}
y_{j+1}=\frac{h}{2}z_{j+1},~
y_{j}=-\frac{h}{2}z_{j},~j=0,1,\cdots,N.
\end{equation}
  This, combining with $y_{0}=0$ and $y_{N+1}=0$, gives
   $y_{j}=0$ ($j=1,2,\cdots,N$) which  is also a contradiction.
   This completes the proof of the lemma.
   \hfill\rule{2 mm}{3mm}

   \vspace{0.2cm}

    The  following lemma  comes from \cite{Liujk1}.
\begin{lemma}\label{lem3.1}
Let $\{u_{i}\}_{i}$, $\{v_{i}\}_{i}$ and $\{w_{i}\}_{i}$ be the sequences of complex numbers. Then,
\begin{equation}
\begin{array}{l}
\disp
\quad\frac{1}{4}\sum_{i=0}^{N}(u_{i+1}-u_{i})(v_{i+1}+v_{i})(w_{i+1}+w_{i})\crr\disp+
\frac{1}{4}\sum_{i=0}^{N}(u_{i+1}-u_{i})(v_{i+1}-v_{i})(w_{i+1}-w_{i})\crr\disp
+\frac{1}{4}\sum_{i=0}^{N}(u_{i+1}+u_{i})(v_{i+1}-v_{i})(w_{i+1}+w_{i})\crr\disp
+\frac{1}{4}\sum_{i=0}^{N}(u_{i+1}+u_{i})(v_{i+1}+v_{i})(w_{i+1}-w_{i})\crr\disp
=u_{N+1}v_{N+1}w_{N+1}-u_{0}v_{0}w_{0}.
\end{array}
\end{equation}
\end{lemma}

The following uniformly stability criterion which was presented  in \cite{Liuzy} or \cite{Abdallah1}
   will be used in  the proof of our main result  Theorem \ref{th2.3} later.

\begin{theorem}\label{th2.2}
 Let $h^{*}>0$ and let $\{S_{h}(t)\}_{h\in (0,h^*)}$ be a family of semigroups of contractions on the Hilbert
space $H_h$, and let $\widetilde{A}_{h}$ be the corresponding infinitesimal generators.
 The family $\{S_{h}(t)\}$  is uniformly exponentially stable if and only if the   following two conditions are fulfilled:

\begin{itemize}

\item For every   $h\in(0,h^{*})$, $i\mathbb{R}\subset \rho(\widetilde{A}_{h})$;

\item $\sup_{h\in(0,h^{*}),\beta\in\mathbb{R}}\|(i\beta I-\widetilde{A}_{h})^{-1}\|<\infty$.
\end{itemize}
\end{theorem}

Now, we are in a  position to give the main result of this paper.

\begin{theorem}\label{th2.3} For the matrices  $\mathcal{A}_{h}$ defined by (\ref{A}),
the corresponding  family of $C_0$-semigroups  $T_{h}(t)$ generated by $\mathcal{A}_{h}$
is uniformly exponentially stable, that is, there exist two constants $M > 0$ and  $ \omega> 0$ independent of $h \in (0, 1)$
 such that
 \begin{equation}\label{3.12}
 \|T_{h}(t)\|\leq M e^{-\omega t},~\forall t\geq0.
 \end{equation}
\end{theorem}

{\bf Proof.} The proof is based on Theorem \ref{th2.2}. Notice that by Lemma \ref{pro3.1+},   for every  $h\in(0,1)$,
  $T_{h}(t)$ is a  $C_0$-semigroup of contractions.
  The fact that $\mathcal{A}_{h}$ satisfies the first condition  of  Theorem \ref{th2.2} has been claimed by Lemma \ref{pro3.2}. In order
to show that the family    $\mathcal{A}_{h}$ satisfies the second condition of  Theorem \ref{th2.2},
 we prove by  contradiction.
 If the second condition  of  Theorem \ref{th2.2} is false, then there exist a sequence  $\beta_{n}\in \mathbb{R}$,
 $h_{n}\in (0,1)$, and  $Y^{n}_{h_{n}}\in \mathbb{Y}_{h_{n}}, \|Y^{n}_{h_{n}}\|_{\mathbb{Y}_{h_{n}}}=1$  such that
\begin{equation}\label{3.13}
\|U^{n}_{h_{n}}\|_{\mathbb{Y}_{h_{n}}}\leq n^{-1}, \; U^{n}_{h_{n}}
=(i\beta_{n} I_{h_{n}}-\mathcal{A}_{h_{n}})Y^{n}_{h_{n}}.
\end{equation}
By  the Cauchy-Schwartz inequality, it follows from (\ref{3.13}) and (\ref{3.5}) that
\begin{equation}\label{3.18}
\begin{array}{ll}
\disp \textrm{Re}\left\langle U^{n}_{h_{n}}, Y^{n}_{h_{n}}\right\rangle_{\mathbb{Y}_{h_{n}}}&=\disp -\textrm{Re}\left\langle \mathcal{A}_{h_{n}}Y^{n}_{h_{n}}, Y^{n}_{h_{n}}\right\rangle_{\mathbb{Y}_{h_{n}}}\crr
&=k|y^{n}_{N_{n}+1}|^{2}\leq n^{-1}.
\end{array}
\end{equation}

Let  $$Z_{h_{n}}^{n}=(z_{0}^{n}, z_{1}^{n},\cdots, z_{N_{n}}^{n})^{\top}\in \mathbb{Z}_{h_{n}}$$ be the shadow element of $Y_{h_{n}}^{n}=(y_{1}^{n}, y_{2}^{n},\cdots, y_{N_{n}+1}^{n})^{\top} {\rm(see \dref{zh1})}$,  $U_{h_{n}}^{n}=(u_{1}^{n}, u_{2}^{n},\cdots, u_{N_{n}+1}^{n})^{\top}$ with $h_{n}(N_{n}+1)=1$.
 Set artificially
 $u_{0}^{n}=y_{0}^{n}=0$ and $z_{N_{n}+1}^{n}=-iky_{N_{n}+1}^{n}$
 to unify  the  notation of  $u^{n}_{j+\frac{1}{2}}$ and $\delta_{x}z^{n}_{j+\frac{1}{2}}$
from $j=0,1,\cdots,N_{n}$. Then, it follows from (\ref{3.13}) that
 \begin{equation}\label{3.14}
 \left\{\begin{array}{l}
 \disp
D_{h_{n}}U^{n}_{h_{n}}
=i\beta_{n}D_{h_{n}}Y^{n}_{h_{n}}+iM_{h_{n}}Z^{n}_{h_{n}}+(0,\cdots,0,ih^{-1}z_{N_{n}+1}^{n})^{\top},\crr\disp
-M_{h_{n}}^{\top}Y^{n}_{h_{n}}=D_{h_{n}}^{\top}Z^{n}_{h_{n}}+(0,\cdots,0,2^{-1}z_{N_{n}+1}^{n})^{\top},
\end{array}\right.
\end{equation}
or in vector form:
 \begin{equation}\label{3.16}
 \left\{\begin{array}{ll}
 \disp
\left(\begin{array}{l}
~u^{n}_{0+\frac{1}{2}}\\
~u^{n}_{1+\frac{1}{2}}\\
~~~\vdots\\
u^{n}_{N_{n}+\frac{1}{2}}
\end{array}
\right)&=\disp
i\beta_{n}\left(\begin{array}{l}
~y^{n}_{0+\frac{1}{2}}\\
~y^{n}_{1+\frac{1}{2}}\\
~~~\vdots\\
y^{n}_{N_{n}+\frac{1}{2}}
\end{array}
\right)+
i\left(\begin{array}{l}
~\delta_{x}z^{n}_{0+\frac{1}{2}}\\
~\delta_{x}z^{n}_{1+\frac{1}{2}}\\
~~~\vdots\\
\delta_{x}z^{n}_{N_{n}+\frac{1}{2}}
\end{array}
\right),\crr\disp
\disp
\left(\begin{array}{l}
~z^{n}_{0+\frac{1}{2}}\\
~z^{n}_{1+\frac{1}{2}}\\
~~~\vdots\\
z^{n}_{N_{n}+\frac{1}{2}}
\end{array}
\right)&=\disp \left(\begin{array}{l}
~\delta_{x}y^{n}_{0+\frac{1}{2}}\\
~\delta_{x}y^{n}_{1+\frac{1}{2}}\\
~~~\vdots\\
\delta_{x}y^{n}_{N_{n}+\frac{1}{2}}
\end{array}
\right).
\end{array}\right.
\end{equation}

 The proof will be split into three claims and each claim corresponds to that in the proof of stability of PDE.
 Clam 1  corresponds to $\omega_{n}\rightarrow\infty$ in the proof of Theorem \ref{Th.A1}.

{\bf Cliam 1:    $|\beta_n|\geq C^\prime>0$ for some constant $C^\prime$ independent of $ n \in  \mathbb{N}^+$.}

  Suppose by contrary that the sequence $\{\beta_{n}\}$ contains a subsequence
which is still denoted by $\{\beta_{n}\}$ itself  without loss of generality
converging  to zero.  Since  $\|Y^{n}_{h_{n}}\|_{\mathbb{Y}_{h_{n}}}=1$ and
$\|U^{n}_{h_{n}}\|_{\mathbb{Y}_{h_{n}}}\leq n^{-1}$, it follows from (\ref{3.16}) that
\begin{equation}\label{3.19}
\begin{array}{l}
 \disp\quad h_{n}\sum_{j=0}^{N_{n}}\left|\delta_{x}z^{n}_{j+\frac{1}{2}}\right|^{2}=
h_{n}\sum_{j=0}^{N_{n}}\left|u^{n}_{j+\frac{1}{2}}-i\beta_{n}y^{n}_{j+\frac{1}{2}}\right|^{2}\crr\disp
\leq  \disp 2h_{n}\sum_{j=0}^{N_{n}}\left|u^{n}_{j+\frac{1}{2}}\right|^{2}+
2\beta_{n}^{2}h_{n}\sum_{j=0}^{N_{n}}\left|y^{n}_{j+\frac{1}{2}}\right|^{2}\crr\disp
=2\|U^{n}_{h_{n}}\|_{\mathbb{Y}_{h_{n}}}^{2}+
2\beta_{n}^{2}\|Y^{n}_{h_{n}}\|_{\mathbb{Y}_{h_{n}}}^{2}\leq 2n^{-2},
\end{array}
\end{equation}
which holds for all sufficiently large $n$.  On the other hand, by some simple operations,   we get
\begin{equation}
\begin{array}{l}
\disp \quad |z_{j}^{n}-z_{N_{n}+1}^{n}|^{2}=|z_{j}^{n}-z_{j+1}^{n}+z_{j+1}^{n}-z_{j+2}^{n}+z_{j+2}^{n}\cdots-z_{N_{n}+1}^{n}|^{2}
\crr\disp
=\left|\sum_{l=j}^{N_{n}}(z_{l+1}^{n}-z_{l}^{n})\right|^{2}
\leq\left(\sum_{l=j}^{N_{n}}\left|1\right|^{2}\right)\left(\sum_{l=j}^{N_{n}}\left|z_{l+1}^{n}-z_{l}^{n}\right|^{2}\right)
\crr\disp\leq (N_{n}+1)\left(\sum_{l=0}^{N_{n}}\left|z_{l+1}^{n}-z_{l}^{n}\right|^{2}\right)\crr\disp
=h_{n}\sum_{j=1}^{N_{n}}\left|\delta_{x}z_{j+\frac{1}{2}}^{n}\right|^{2},~~j=0,1,\cdots,N_{n},
\end{array}
\end{equation}
in which $ h_{n}(N_{n}+1)=1$ is used in the last step, and for $j=0,1,\cdots,N_{n}$
$$
|z_{j}^{n}|\leq|z_{j}^{n}-z_{N_{n}+1}^{n}|+|z_{N_{n}+1}^{n}|\leq\sqrt{h_{n}
\sum_{j=1}^{N_{n}}|\delta_{x}z_{j+\frac{1}{2}}^{n}|^{2}}+|z_{N_{n}+1}^{n}|.
$$
This inequality, together with $z_{N_{n}+1}^{n}=-iky_{N_{n}+1}^{n}$ and (\ref{3.18})-(\ref{3.19}),
 implies that for each $j=0,1,\cdots,N_{n}$, $|z^{n}_{j}|^{2}=\mathcal{O}(n^{-1})$.
  Therefore, in light of $h_{n}(N_{n}+1)=1$ and the second identity of (\ref{3.16}),
\begin{equation}\label{3.23}
\begin{array}{l}
\disp\quad h_{n}\sum_{j=0}^{N_{n}}
\left|z^{n}_{j+\frac{1}{2}}\right|^{2}\leq \frac{h_{n}}{2}\sum_{j=0}^{N_{n}}\left(
\left|z^{n}_{j+1}\right|^{2}+\left|z^{n}_{j}\right|^{2}\right)\crr
\disp
\leq h_{n}(N_{n}+1)\left|\sqrt{h_{n}
\sum_{j=1}^{N_{n}}|\delta_{x}z_{j+\frac{1}{2}}^{n}|^{2}}+|z_{N_{n}+1}^{n}|\right|^{2}
\crr
\disp\leq h_{n}(N_{n}+1)C n^{-1}=\mathcal{O}(n^{-1}).
\end{array}
\end{equation}
Thus, the deducing process from (\ref{3.19}) to (\ref{3.23}) tells us that
$$
h_{n}\sum_{j=0}^{N_{n}}
\left|\delta_{x}z^{n}_{j+\frac{1}{2}}\right|^{2}=\mathcal{O}(n^{-2}),
$$
which implies that
$$
h_{n}\sum_{j=0}^{N_{n}}
\left|z^{n}_{j+\frac{1}{2}}\right|^{2}=\mathcal{O}(n^{-1}).
$$

By noticing  the second identity of (\ref{3.16}), we have
$$
h_{n}\sum_{j=0}^{N_{n}}
\left|\delta_{x}y^{n}_{j+\frac{1}{2}}\right|^{2}=h_{n}\sum_{j=0}^{N_{n}}
\left|z^{n}_{j+\frac{1}{2}}\right|^{2},
$$
which means, by (\ref{3.23})  that
$$
h_{n}\sum_{j=0}^{N_{n}}
\left|\delta_{x}y^{n}_{j+\frac{1}{2}}\right|^{2}=\mathcal{O}(n^{-1}).
$$
Similarly, repeating the procedure from (\ref{3.19}) to (\ref{3.23}),  for $Y_{h_{n}}^{n}$,
 we obtain
 $$
 \|Y^{n}_{h_{n}}\|_{\mathbb{Y}_{h_{n}}}^{2}=h_{n}\sum_{j=0}^{N_{n}}
\left|y^{n}_{j+\frac{1}{2}}\right|^{2}=\mathcal{O}(n^{-1/2}),
$$
 which leads to a contradiction.  Thus, the sequence $\{\beta_{n}\}$ cannot contain a subsequence converging  to zero.
 We can therefore assume that  $|\beta_{n}|\geq C'>0$ for some constant $C'$ independent of $n\in \mathbb{N}^{+}$.

The second  claim is the discrete counterpart of (\ref{1.7+}) but with two extra terms
$$
\frac{h_{n}^{3}}{4\beta_{n}}\sum_{j=0}^{N_{n}}
\left|\delta_{x}z^{n}_{j+\frac{1}{2}}\right|^{2}+
\frac{h_{n}^{3}}{4}\sum_{j=0}^{N_{n}}\left|\delta_{x}y^{n}_{j+\frac{1}{2}}\right|^{2}
$$
which play important  roles in our proofs.

{\bf Claim 2: the following \dref{gbz1} holds true:
\begin{eqnarray}\label{gbz1}
&&\disp \quad \quad\|Y^{n}_{h_{n}}\|_{\mathbb{Y}_{h_{n}}}^{2}+\frac{1}{\beta_{n}}h_n\|\Sigma_{h_n}\widehat{Z}^{n}_{h_{n}}\|_{\mathbb{C}^{N_{n}+2}}^2
\crr&&\disp\quad+\frac{h_{n}^{2}}{4\beta_{n}}h_n\|\Delta_{h_n}\widehat{Z}^{n}_{h_{n}}\|_{\mathbb{C}^{N_{n}+2}}^{2}
+\frac{h_{n}^{2}}{4}h_n\|\Delta_{h_n}\widehat{Y}^{n}_{h_{n}}\|_{\mathbb{C}^{N_{n}+2}}^{2} \crr
&&\disp
=h_{n}\sum_{j=0}^{N_{n}}
\left|y^{n}_{j+\frac{1}{2}}\right|^{2}+\frac{h_{n}}{\beta_{n}}\sum_{j=0}^{N_{n}}
\left|z^{n}_{j+\frac{1}{2}}\right|^{2}\crr
&&\disp\quad+\frac{h_{n}^{3}}{4\beta_{n}}\sum_{j=0}^{N_{n}}
\left|\delta_{x}z^{n}_{j+\frac{1}{2}}\right|^{2}+
\frac{h_{n}^{3}}{4}\sum_{j=0}^{N_{n}}\left|\delta_{x}y^{n}_{j+\frac{1}{2}}\right|^{2}
\crr&&\disp=\mathcal{O}(n^{-1}),
\end{eqnarray}}
where $\|\cdot\|_{\mathbb{C}^{N_{n}+2}}$ denotes the standard norm
of $\mathbb{C}^{N_{n}+2}$ and
\begin{equation}\label{dm2}
\begin{array}{l}
\disp
\widehat{Z}^{n}_{h_{n}}=(z_{0},z_{1},\cdots,z_{N_{n}+1})^{\top}=\left((Z^{n}_{h_{n}})^{\top},z_{N_{n}+1}\right)^{\top}, \crr
\disp
\widehat{Y}^{n}_{h_{n}}=(0,y_{1},\cdots,y_{N_{n}+1})^{\top}=\left(0,(Y^{n}_{h_{n}})^{\top}\right)^{\top},\crr\disp
\Sigma_{h}=\frac{1}{2}\left(\begin{matrix}
1&1&&&\\
&\ddots&\ddots&&\\
&&1&1&\\
&&&1&1
\end{matrix}\right)_{(N+2)\times(N+1)},\crr\disp
\Delta_{h}=\frac{1}{h}\left(\begin{matrix}
-1&1&&&\\
&-1&1&&\\
&&\ddots&\ddots&\\
&&&-1&1\\
\end{matrix}\right)_{(N+2)\times(N+1)}.
\end{array}
\end{equation}

Actually, it follows from (\ref{3.13}),   (\ref{3.16}), and  $\|Y^{n}_{h_{n}}\|_{\mathbb{Y}_{h_{n}}}=1$ that $\beta_{n}^{-2}h_{n}\sum_{j=0}^{N_{n}}
\left|\delta_{x}z^{n}_{j+\frac{1}{2}}\right|^{2}$ is uniformly bounded with respect to $n\in\mathbb{N}^{+}$ because
\begin{align}
\frac{1}{\beta_{n}}\left(\begin{array}{l}
~\delta_{x}z^{n}_{0+\frac{1}{2}}\\
~\delta_{x}z^{n}_{1+\frac{1}{2}}\\
~~~\vdots\\
\delta_{x}z^{n}_{N_{n}+\frac{1}{2}}
\end{array}
\right)=-\left(\begin{array}{l}
~y^{n}_{0+\frac{1}{2}}\\
~y^{n}_{1+\frac{1}{2}}\\
~~~\vdots\\
y^{n}_{N_{n}+\frac{1}{2}}
\end{array}
\right)-\frac{i}{\beta_{n}}\left(\begin{array}{l}
~u^{n}_{0+\frac{1}{2}}\\
~u^{n}_{1+\frac{1}{2}}\\
~~~\vdots\\
u^{n}_{N_{n}+\frac{1}{2}}
\end{array}
\right).\nonumber
\end{align}
By (\ref{3.23})  and {\bf Claim  1},
$$
\beta_{n}^{-2}h_{n}\sum_{j=0}^{N_{n}}
\left|z^{n}_{j+\frac{1}{2}}\right|^{2}
$$
 is also uniformly bounded with respect to $n\in\mathbb{N}^{+}$.
  Let $x_{j}^{n}=jh_{n}$ for $j=0,1,\cdots,N_{n}+1$ and consider the following estimates:
\begin{equation}\label{3.24}
\begin{array}{ll}
&\disp \quad\left|h_{n}\sum_{j=0}^{N_{n}}x_{j+\frac{1}{2}}^{n}(\beta_{n}y^{n}_{j+\frac{1}{2}}+\delta_{x}z^{n}_{j+\frac{1}{2}})
\frac{\overline{z}^{n}_{j+\frac{1}{2}}}{\beta_{n}}\right|^{2}\crr\disp
&=\disp
\left|h_{n}\sum_{j=0}^{N_{n}}x_{j+\frac{1}{2}}^{n}u^{n}_{j+\frac{1}{2}}
\frac{\overline{z}^{n}_{j+\frac{1}{2}}}{\beta_{n}}\right|^{2}\crr\disp
& \leq\disp
\left(\sum_{j=0}^{N_{n}}\left|\sqrt{h_{n}}u^{n}_{j+\frac{1}{2}}\right|
\left|\sqrt{h_{n}}\frac{z^{n}_{j+\frac{1}{2}}}{\beta_{n}}\right|\right)^{2}\crr
&\disp  \leq \left(h_{n}\sum_{j=0}^{N_{n}}
\left|u^{n}_{j+\frac{1}{2}}\right|^{2}\right)\left(\beta_{n}^{-2}h_{n}\sum_{j=0}^{N_{n}}
\left|z^{n}_{j+\frac{1}{2}}\right|^{2}\right)\crr
&\disp=\|U^{n}_{h_{n}}\|_{\mathbb{Y}_{h_{n}}}^{2}
\left(\beta_{n}^{-2}h_{n}\sum_{j=0}^{N_{n}}
\left|z^{n}_{j+\frac{1}{2}}\right|^{2}\right)\crr
&\disp=\mathcal{O}(n^{-2}),
\end{array}
\end{equation}
where  (\ref{3.13}) and (\ref{3.16}) were used. On the other hand, by the second identity of (\ref{3.16}), we have
\begin{equation}\label{3.25}
\begin{array}{l}
\disp \quad h_{n}\sum_{j=0}^{N_{n}}x_{j+\frac{1}{2}}^{n}(\beta_{n}y^{n}_{j+\frac{1}{2}}+\delta_{x}z^{n}_{j+\frac{1}{2}})
\frac{\overline{z}^{n}_{j+\frac{1}{2}}}{\beta_{n}}\crr\disp
\disp =h_{n}\sum_{j=0}^{N_{n}}x_{j+\frac{1}{2}}^{n}y^{n}_{j+\frac{1}{2}}\delta_{x}\overline{y}^{n}_{j+\frac{1}{2}}\crr\disp
\disp \quad+\beta_{n}^{-1}h_{n}\sum_{j=0}^{N_{n}}x_{j+\frac{1}{2}}^{n}\delta_{x}z^{n}_{j+\frac{1}{2}}\overline{z}^{n}_{j+\frac{1}{2}}.
\end{array}
\end{equation}
Applying Lemma \ref{lem3.1} to the two terms of the right hand side of \dref{3.25} and noticing  $x_{N_{n}+1}^{n}=1$, $x_{0}^{n}=0$, $x_{j+1}^{n}-x_{j}^{n}=h_{n}$, it is easy to obtain
\begin{equation}\label{3.26}
\begin{array}{l}
\disp
\quad 2{\rm Re}\left(h_{n}\sum_{j=0}^{N_{n}}x_{j+\frac{1}{2}}^{n}y^{n}_{j+\frac{1}{2}}\delta_{x}\overline{y}^{n}_{j+\frac{1}{2}}\right)
\crr
=\disp h_{n}\sum_{j=0}^{N_{n}}x_{j+\frac{1}{2}}^{n}y^{n}_{j+\frac{1}{2}}\delta_{x}\overline{y}^{n}_{j+\frac{1}{2}}
+h_{n}\sum_{j=0}^{N_{n}}x_{j+\frac{1}{2}}^{n}\overline{y}^{n}_{j+\frac{1}{2}}\delta_{x}y^{n}_{j+\frac{1}{2}} \crr\disp= \frac{1}{4}\sum_{j=0}^{N_{n}}(x_{j+1}^{n}+x_{j}^{n})(y^{n}_{j+1}+y^{n}_{j})(\overline{y}^{n}_{j+1}-\overline{y}^{n}_{j})
\crr\disp \quad+\frac{1}{4}\sum_{j=0}^{N_{n}}(x_{j+1}^{n}+x_{j}^{n})(y^{n}_{j+1}-y^{n}_{j})(\overline{y}^{n}_{j+1}+\overline{y}^{n}_{j})\crr
\disp  \disp = \left|y_{N_{n}+1}\right|^{2}-h_{n}\sum_{j=0}^{N_{n}}\left|y^{n}_{j+\frac{1}{2}}\right|^{2}
-\frac{h_{n}^{3}}{4}\sum_{j=0}^{N_{n}}\left|\delta_{x}y^{n}_{j+\frac{1}{2}}\right|^{2},
\end{array}
\end{equation}
and
\begin{equation}\label{3.27}
\begin{array}{l}
\disp \quad 2{\rm Re}\left(h_{n}\sum_{j=0}^{N_{n}}x_{j+\frac{1}{2}}^{n}z^{n}_{j+\frac{1}{2}}\delta_{x}\overline{z}^{n}_{j+\frac{1}{2}}\right)
\crr=\disp
h_{n}\sum_{j=0}^{N_{n}}x_{j+\frac{1}{2}}^{n}z^{n}_{j+\frac{1}{2}}\delta_{x}\overline{z}^{n}_{j+\frac{1}{2}}
+h_{n}\sum_{j=0}^{N_{n}}x_{j+\frac{1}{2}}^{n}\overline{z}^{n}_{j+\frac{1}{2}}\delta_{x}z^{n}_{j+\frac{1}{2}}
\crr\disp=\frac{1}{4}\sum_{j=0}^{N_{n}}(x_{j+1}^{n}+x_{j}^{n})(z^{n}_{j+1}+z^{n}_{j})(\overline{z}^{n}_{j+1}-\overline{z}^{n}_{j})\crr
\disp~~+\frac{1}{4}\sum_{j=0}^{N_{n}}(x_{j+1}^{n}+x_{j}^{n})(\overline{z}^{n}_{j+1}+\overline{z}^{n}_{j})(z^{n}_{j+1}-z^{n}_{j})\crr
=\disp |z_{N_{n}+1}|^{2}-h_{n}\sum_{j=0}^{N_{n}}\left|z^{n}_{j+\frac{1}{2}}\right|^{2}
-\frac{h_{n}^{3}}{4}\sum_{j=0}^{N_{n}}\left|\delta_{x}z^{n}_{j+\frac{1}{2}}\right|^{2}.
\end{array}
\end{equation}
By  (\ref{3.25})-(\ref{3.27}), it follows that
\begin{equation}
\begin{array}{ll}
&\disp
 \quad h_{n}\sum_{j=0}^{N_{n}}
\left|y^{n}_{j+\frac{1}{2}}\right|^{2}+
\frac{h_{n}^{3}}{4}\sum_{j=0}^{N_{n}}\left|\delta_{x}y^{n}_{j+\frac{1}{2}}\right|^{2}\crr
&\disp\quad+\frac{h_{n}}{\beta_{n}}\sum_{j=0}^{N_{n}}
\left|z^{n}_{j+\frac{1}{2}}\right|^{2}+\frac{h_{n}^{3}}{4\beta_{n}}\sum_{j=0}^{N_{n}}
\left|\delta_{x}z^{n}_{j+\frac{1}{2}}\right|^{2}\crr
&\disp
=-2{\rm{Re}}\left(h_{n}\sum_{j=0}^{N_{n}}x_{j+\frac{1}{2}}^{n}(\beta_{n}y^{n}_{j+\frac{1}{2}}+\delta_{x}z^{n}_{j+\frac{1}{2}})
\frac{\overline{z}^{n}_{j+\frac{1}{2}}}{\beta_{n}}\right)\crr
&\disp\quad+|y_{N_{n}+1}|^{2}+\beta_{n}^{-1}|z_{N_{n}+1}|^{2},
\end{array}
\end{equation}
which proves  \dref{gbz1} by  (\ref{3.18}),   (\ref{3.24}) and $z_{N_{n}+1}=-iky_{N_{n}+1}$.

 The third claim  is perfectly the discrete counterpart of (\ref{1.7++}).

{\bf Claim 3: The following \dref{gbz2} holds true:
\begin{equation}\label{gbz2}
\begin{array}{l}
\disp
\quad \|Y^{n}_{h_{n}}\|_{\mathbb{Y}_{h_{n}}}^{2}
+\frac{1}{\beta_{n}^{2}}h_n\|\Delta_{h_n}\widehat{Z}^{n}_{h_{n}}\|_{\mathbb{C}^{N_{n}+2}}^{2}
-\frac{2}{\beta_{n}}h_n\|\Sigma_{h_n}\widehat{Z}^{n}_{h_{n}}\|_{\mathbb{C}^{N_{n}+2}}^{2}
\crr\disp=h_{n}\sum_{j=0}^{N_{n}}
\left|y^{n}_{j+\frac{1}{2}}\right|^{2}+\frac{h_{n}}{\beta_{n}^{2}}\sum_{j=0}^{N_{n}}
\left|\delta_{x}z^{n}_{j+\frac{1}{2}}\right|^{2}-\frac{2h_{n}}{\beta_{n}}\sum_{j=0}^{N_{n}}
\left|z^{n}_{j+\frac{1}{2}}\right|^{2}
\crr
\disp =\mathcal{O}(n^{-2}),
\end{array}
\end{equation}
in which  $\Sigma_{h_n}$  and  $\Delta_{h_n}$  are defined  in \dref{dm2} and $\|\cdot\|_{\mathbb{C}^{N_{n}+2}}$ denotes the standard norm
of $\mathbb{C}^{N_{n}+2}$. }

Actually, from (\ref{3.16}), we have
\begin{equation}\label{3.29}
\begin{array}{l}
\disp \quad \frac{\|U^{n}_{h_{n}}\|_{\mathbb{Y}_{h_{n}}}^{2}}{\beta_{n}^{2}}\crr=\disp \frac{h_{n}}{\beta_{n}^{2}}\sum_{j=0}^{N_{n}}
\left|u^{n}_{j+\frac{1}{2}}\right|^{2}=h_{n}\sum_{j=0}^{N_{n}}
\left|y^{n}_{j+\frac{1}{2}}+\frac{\delta_{x}z^{n}_{j+\frac{1}{2}}}{\beta_{n}}\right|^{2}\crr
=\disp h_{n}\sum_{j=0}^{N_{n}}
\left|y^{n}_{j+\frac{1}{2}}\right|^{2}+\frac{h_{n}}{\beta_{n}^{2}}\sum_{j=0}^{N_{n}}
\left|\delta_{x}z^{n}_{j+\frac{1}{2}}\right|^{2}\crr\disp~~~+\frac{h_{n}}{\beta_{n}}\sum_{j=0}^{N_{n}}(
\overline{y}^{n}_{j+\frac{1}{2}}\delta_{x}z^{n}_{j+\frac{1}{2}}
+y^{n}_{j+\frac{1}{2}}\delta_{x}\overline{z}^{n}_{j+\frac{1}{2}}).
\end{array}
\end{equation}
On the other hand, it follows from the second identity of (\ref{3.16}) that $z^{n}_{j+\frac{1}{2}}=\delta_{x}y^{n}_{j+\frac{1}{2}}$ and
\begin{eqnarray*}
&&\disp \quad\frac{h_{n}}{\beta_{n}}\sum_{j=0}^{N_{n}}(
\overline{y}^{n}_{j+\frac{1}{2}}\delta_{x}z^{n}_{j+\frac{1}{2}}
+y^{n}_{j+\frac{1}{2}}\delta_{x}\overline{z}^{n}_{j+\frac{1}{2}})+
\frac{2h_{n}}{\beta_{n}}\sum_{j=0}^{N_{n}}\left|z^{n}_{j+\frac{1}{2}}\right|^{2}\crr
&&\disp =\frac{h_{n}}{\beta_{n}}\sum_{j=0}^{N_{n}}(
\overline{y}^{n}_{j+\frac{1}{2}}\delta_{x}z^{n}_{j+\frac{1}{2}}
+\delta_{x}\overline{y}^{n}_{j+\frac{1}{2}}z^{n}_{j+\frac{1}{2}})
\crr
&&\disp\quad+\frac{h_{n}}{\beta_{n}}\sum_{j=0}^{N_{n}}
(\delta_{x}y^{n}_{j+\frac{1}{2}}\overline{z}^{n}_{j+\frac{1}{2}}
+y^{n}_{j+\frac{1}{2}}\delta_{x}\overline{z}^{n}_{j+\frac{1}{2}})\crr
&&\disp =\frac{1}{2\beta_{n}}\sum_{j=0}^{N_{n}}[(
\overline{y}^{n}_{j+1}+\overline{y}^{n}_{j})(z^{n}_{j+1}-z^{n}_{j})+(
\overline{y}^{n}_{j+1}-\overline{y}^{n}_{j})(z^{n}_{j+1}+z^{n}_{j})]\crr
&&\disp \quad+\frac{1}{2\beta_{n}}
\sum_{j=0}^{N_{n}}[(y^{n}_{j+1}+y^{n}_{j})(\overline{z}^{n}_{j+1}-\overline{z}^{n}_{j})+
(y^{n}_{j+1}-y^{n}_{j})(\overline{z}^{n}_{j+1}+\overline{z}^{n}_{j})]
\crr
&&\disp =\frac{1}{\beta_{n}}\sum_{j=0}^{N_{n}}(
\overline{y}^{n}_{j+1}z^{n}_{j+1}-\overline{y}^{n}_{j}z^{n}_{j})+\frac{1}{\beta_{n}}
\sum_{j=0}^{N_{n}}(y^{n}_{j+1}\overline{z}^{n}_{j+1}-y^{n}_{j}\overline{z}^{n}_{j})
\crr
&&\disp =\frac{1}{\beta_{n}}[\overline{y}^{n}_{N_{n}+1}z^{n}_{N_{n}+1}+y^{n}_{N_{n}+1}\overline{z}^{n}_{N_{n}+1}
-\overline{y}^{n}_{0}z^{n}_{0}-y^{n}_{0}\overline{z}^{n}_{0}]=0,
\end{eqnarray*}
 where $z^{n}_{N_{n}+1}=-iky^{n}_{N_{n}+1}$ and $y_{0}^{n}=0$ were used  in the last step. Hence
  \begin{equation}\label{pp1}
\frac{h_{n}}{\beta_{n}}\sum_{j=0}^{N_{n}}(
\overline{y}^{n}_{j+\frac{1}{2}}\delta_{x}z^{n}_{j+\frac{1}{2}}
+y^{n}_{j+\frac{1}{2}}\delta_{x}\overline{z}^{n}_{j+\frac{1}{2}})
=-\frac{2h_{n}}{\beta_{n}}\sum_{j=0}^{N_{n}}\left|z^{n}_{j+\frac{1}{2}}\right|^{2}.
\end{equation}
   Plugging \dref{pp1}  into (\ref{3.29}) and using  $\|U^{n}_{h_{n}}\|_{\mathbb{Y}_{h_{n}}}\leq n^{-1}$ and
   {\bf Claim 1}, we arrive at \dref{gbz2}.

   Finally,  if $\beta_{n}>0$,   we have $\|Y^{n}_{h_{n}}\|_{\mathbb{Y}_{h_{n}}}^{2}=\mathcal{O}(n^{-1})$ from  \dref{gbz1},
 which is  a contradiction to $\|Y^{n}_{h_{n}}\|_{\mathbb{Y}_{h_{n}}}=1$.
  When $\beta_{n}<0$,   $\|Y^{n}_{h_{n}}\|_{\mathbb{Y}_{h_{n}}}=\mathcal{O}(n^{-1})$  by virtue of \dref{gbz2},
  which is  also a contradiction. We therefore complete the
  proof of the theorem.   \hfill\rule{2 mm}{3mm}

\section{Concluding remarks}\label{sec5}

In this paper, the uniform approximation of exponential stability
of a one-dimensional Schr\"{o}dinger equation is investigated.
We introduce an order reduction  space semi-discretized finite difference scheme
 for approximating uniformly the exponentially stable closed-loop system.
 Although the scheme has been applied to certain PDEs in previous works, they all share a common feature
that  it is possible to find a suitable Lyapunov functional for the closed-loop systems, both for the continuous system and its discrete counterpart.
 However, for the system considered in this paper,  it  is a longstanding  problem
 that  in the natural state space  $L^2(0,1)$, even for the continuous system, the time domain energy multiplier
 has not been found.
 This makes the convergence of the semi-discrete scheme of this PDE  be open for  a long time.
 This paper is the first work that applies the frequency domain multiplier approach to the uniformly exponential convergence of semi-discretized  PDE system.
The convergence of the discrete scheme to continuous
 system is not included because it is a standard procedure and can be followed analogously from \cite{Liujk2} and many other papers.
 Considering it is  difficult to find a time domain  energy multiplier for many other PDEs, the approach presented in this paper has  significant potentials in applying to other PDEs.

\section*{Acknowledgments}
The authors would like to thank Dr. Jiankang Liu and Dr. Hanjin Ren for their careful
reading and many suggestions on the presentation of the
initial manuscript of the paper.  The Figures 1 and 2
were depicted by Dr. Jiangkang Liu.
\ifCLASSOPTIONcaptionsoff
  \newpage
\fi





\begin{thebibliography}{1}




\bibitem{Abdallah1} F. Abdallah, S. Nicaise, J. Valein, and A. Wehbe, Uniformly exponentially or polynomially stable
approximations for second order evolution equations and some applications,
 {\it ESAIM Control Optim. Calc. Var.}, 19(2013), 844-887.

 \bibitem{Banks}H.T. Banks, K. Ito, and C. Wang,
 Exponentially stable approximations of weakly damped wave equations,
  in: \emph{``Estimation and Control of Distributed Parameter Systems''}, Basel, Birkhauser, 1991, pp.1-33.


\bibitem{Dai}L. Dai,  \emph{Singular Control Systems}, SpringerVerlag, New York, 1989.

 \bibitem{Ervedoza1}S. Ervedoza, A. Marica, and E. Zuazua, Numerical meshes
 ensuring uniform observability of 1d waves construction and analysis,
 {\it  IMA J. Numer. Anal.}, 36(2015), 503-542.


\bibitem{Ervedoza2}S. Ervedoza, Spectral conditions for admissibility and
observability of wave systems: applications to finite element schemes,
{\it Numer. Math}, 113(2009), 377-415.

     \bibitem{Ervedoza3}S. Ervedoza, Spectral conditions for admissibility and
observability of Schr\"{o}dinger systems: Applications to finite element discretizations,
{\it  Asymptot. Anal.}, 71(2011), 1-32.

 \bibitem{Glowinski}R. Glowinski, C.H. Li, and J.L. Lions, A numerical approach to the exact
boundary controllability of the wave equation. I. Dirichlet controls: description
of the numerical methods, {\it  Japan J. Appl. Math.}, 7(1990), 1-76.
%

\bibitem{Guobz2}B.Z. Guo and B.B. Xu, A Semi-discrete Finite Difference Method to Uniform Stabilization of Wave Equation with Local Viscosity,
 {\it  IFAC J. Syst. Control }, 13(2020), Art.101000.

\bibitem{Huang}F.L. Huang, Characteristic conditions for exponential stability of linear dynamical
systems in Hilbert spaces, {\it Ann. Differential Equations}, 1(1)(1985), 43-56.

\bibitem{Infante}J.A. Infante and E. Zuazua, Boundary observability for the space semi-discretizations of the 1-d wave equation, \emph{ M2AN Math. Model. Numer. Anal.},  33(1999), 407-438.
%

    \bibitem{Kao}Y. Kao, J. Xie and C. Wang, Stabilization of Singular Markovian Jump Systems With Generally Uncertain Transition Rates, {\it   IEEE Trans. Automat. Control}, 59(2014), 2604-2610.

     \bibitem{Kristic2}M. Krstic, B.Z. Guo,  and A. Smyshlyaev,  Boundary controllers and observers for the linearized Schr\"{o}dinger equation,  SIAM J. Control Optim., 49(2011),  1479-1497.

   \bibitem{Liujk1}J. Liu and B.Z. Guo,
   A new semi-discretized order reduction finite difference scheme for uniform approximation of 1-D wave equation,
    {\it SIAM J. Control Optim.},  58(2020), 2256-2287.

   \bibitem{Liujk2}J. Liu,  R.Q. Hao, and B.Z. Guo,
   Order reduction-based uniform approximation of exponential stability for one-dimensional Schr\"{o}dinger equation,
     {\it Systems Control Lett.}, 160(2022),  Art.105136.





 \bibitem{Liuzy}Z.Y. Liu and  S.M. Zheng, Liu,
  Exponential stability of the semigroup associated with a thermoelastic system,
  {\it Quart. Appl. Math.},   51(1993), 535-545.


\bibitem{Liuzys}Z.Y. Liu and  S.M. Zheng, Uniform exponnential stability and
approximation in control of a thermoelastic system, {\it SIAM J. Control Optim.}, 32(1994), 1226-1246.


%


 \bibitem{Micu}S. Micu and  C. Castro, Boundary controllability of a linear semi-discrete 1-D wave equation derived from a mixed finite element method,  {\it Numer. Math},   102(2006), 413-462.


 \bibitem{Pruss}J. Pr\"{u}ss, On the spectrum of $C_{0}$-semigroups,
 {\it Trans. Amer. Math. Soc.},  284(1984),  847-857.


\bibitem{Tebou1}L.T. Tebou and  E. ZuaZua,
Uniform exponential long time decay for the space semi-discretization of
a locally damped wave equation via an artificial numerical viscosity,
{\it Numer. Math}, 95(2003), 563-598.

\bibitem{Tebou2}L.T. Tebou and  E. Zuazua,
Uniform boundary stabilization of the finite difference space discretization of the 1-d wave equation,
 {\it  Adv. Comput. Math}, 26(2007), 337-365.

    \bibitem{Tucsnak}M. Tucsnak and G. Weiss,
    \emph{Observation and Control for Operator Semigroups}, Birkhauser Verlag, Basel, 2009.

\bibitem{Wang}Y. Wang, Y. Xia, H. Shen and P. Zhou, SMC Design for Robust Stabilization of Nonlinear Markovian Jump Singular Systems,  {\it
 IEEE Trans. Automat. Control},  63(1)(2018),  219-224.



    \bibitem{Zheng}F. Zheng and H. Zhou, State reconstruction of the wave equation with general viscosity and non-collocated observation and control, {\it  J. Math. Anal. Appl.},  502(2020), Art. 125257.



    \bibitem{Zuazua} E. Zuazua,  Propagation, observation, and control of waves approximated by finite difference methods,    {\it SIAM Rev.}, 47(2005), 197-243.

\end{thebibliography}
%

%
\vspace{-15mm}
\begin{IEEEbiography}
	{Bao-Zhu Guo}
	received the Ph.D. degree from the Chinese
	University of Hong Kong in Applied Mathematics in 1991.
	  Since 2000, he
	has been with the Academy of Mathematics and Systems
	Science, the Chinese Academy of Sciences, where he is
	a Research Professor in mathematical system theory.
	From 2004-2019, he was a chair professor in School of
	Computer Science and Applied Mathematics, University of
	the Witwatersrand, South Africa.   He is the author or co-author of the
	books:  ``Stability
	and Stabilization of Infinite Dimensional Systems with Applications''
	(Springer-Verlag, 1999);  ``Active Disturbance Rejection Control for
	Nonlinear Systems: An Introduction'' (John Wiley \& Sons, 2016);  and
	``Control of Wave and Beam PDEs-The Riesz Basis Approach'' (Springer-Verlag, 2019).
	His
	research interests include theory of  infinite-dimensional systems and
	active disturbance rejection control.
	
	Dr. Guo received the One Hundred Talent Program from the Chinese
	Academy of Sciences (1999), and the National Science Fund for Distinguished
	Young Scholars (2003).
	
\end{IEEEbiography}

\vspace{-1mm}
\begin{IEEEbiography}
	{Fu Zheng}
received the B.Sc. and M.Sc. degrees in mathematics from the Yanbian University, Jilin, China in 2002 and 2005, and the Ph.D.
 from Beijing Institute of Information and Control, Beijing, China in 2011.  He
held a postdoctoral research position in  Academy of Mathematics and Systems Science, Academia Sinica,
Beijing, China, from  2013 to 2015.  From 2012-2021,
he was an  associate professor  at Bohai University, Jinzhou, China.
At the end of 2021, he joined the School of
Science, Hainan University, Hainan China, where he is currently a professor.
 His research interests include control theory of infinite-dimensional systems
 and numerical solutions to control problems of distributed parameter systems.
\end{IEEEbiography}








\end{document}